\documentclass{amsart}

\usepackage {epsfig}
\usepackage {graphicx}

 \newtheorem{thm}{Theorem}
 
 \newtheorem{lemma}[thm]{Lemma}
 
 \newtheorem{theorem}[thm]{Theorem}
 \theoremstyle{definition}
 
 \theoremstyle{remark}

 \newcommand{\Real}{\mathbb{R}}

\usepackage{graphicx}
\begin{document}

\title[Sturm-Hurwitz theorem]
{From Rolle's theorem to the Sturm-Hurwitz theorem}

\author{Guy Katriel}

\address{}

\email{haggaik@wowmail.com}

%\subjclass{}

%\keywords{algebraic differential equation, differentially
%algebraic functions. Mathematics Subject Classification: 34A09
%(12H05)}
%\date{}

%\dedicatory{}

%%% ----------------------------------------------------------------------

\begin{abstract}
We present a proof of the Sturm-Hurwitz theorem, using basic
calculus.

\end{abstract}

%%% ----------------------------------------------------------------------
\maketitle
%%% ----------------------------------------------------------------------

Let $f:\Real\rightarrow\Real$ be a continuous $2\pi$-periodic
function, and let us assume that its Fourier series is given by
\begin{equation}
\label{fourier} f(x)\sim \sum_{k=n}^{\infty}{[a_k \cos(kx)+b_k
\sin(kx)]},
\end{equation}
with $a_n^2+b_n^2\neq 0$. In other words we assume that all the
harmonic components of order less than $n$ vanish. The
Sturm-Hurwitz theorem states that

\begin{theorem}
Assuming (\ref{fourier}), the function $f$ has {\it{at least}}
$2n$ distinct zeroes in the interval $[0,2\pi)$.
\end{theorem}
Sturm stated this result for the case of trigonometric
polynomials, while Hurwitz generalized it to periodic functions.
The interest and importance of the Sturm-Hurwitz theorem has been
highlighted in several recent works of V.I. Arnold
\cite{Arnold2,Arnold1}, both in that it is a simple manifestation
of the `topological economy principle', and in that it can be
applied to prove various geometrical results. Interestingly,
Arnold comments that ``There are many known proofs of this Sturm
theorem but all of them are incomprehensible. Of course I can
reproduce them, but you get no intuition from those proofs"
\cite{Arnold2}. Approaches used to prove the Sturm-Hurwitz theorem
of which I am aware are: (i) proof by contradiction based on
orthogonality (\cite{Hurwitz}, \cite{ps} II-141), (ii) proof using
the argument principle of complex analysis (\cite{ps} III-184),
(iii) proof based on the heat equation \cite{Sturm, Polya}, and,
most recently (iv) a `geometrical' proof based on the theory of
hedgehogs, for the case that $f$ is $C^2$ \cite{Maure}.

This note presents yet a proof of the Sturm-Hurwitz theorem, which
is very elementary in that it uses only simple results of calculus
and basic facts on Fourier series (if one restricts to
trigonometric polynomials the proof becomes even simpler, with no
need for any prerequisites on Fourier series). Although I arrived
at this proof independently, I have learned from correspondence
with colleagues subsequent to the publication of the first version
of this note on the Arxiv that it is actually not new. Professor
A. Chambert-Loir has informed me that this proof was used some
10-15 years ago in France as a standard (but difficult) exercise
for second-year undergraduates. Professor Y. Martinez-Maure has
explained to me that the proof presented below is in fact the same
as the proof presented in \cite{Maure}, although this fact is
somewhat disguised by the geometrical tools and terminology -
indeed the aim of \cite{Maure} was precisely to offer a geometric
interpretation of the Sturm-Hurwitz theorem. Thus I make no claim
to originality in this note, but hope that it will be of interest
to teachers and students of calculus as a nice application of the
standard theorems.

Some remarks on variations of the Sturm-Hurwitz theorem for
non-periodic functions will be made at the end of this note, after
presenting the proof.

\vspace{0.5cm}

Rolle's theorem, which tells us that between any two zeroes of a
differentiable function $f:\Real\rightarrow \Real$ there is a zero
of $f'$, plays a key role in our proof of the Sturm-Hurwitz
theorem. Rolle's theorem implies that if $f$ has at least $m$
zeroes, $f'$ has at least $m-1$ zeroes. But for a $2\pi$-periodic
function one can easily see that we have the following stronger
result: if $f$ has at least $m$ zeroes in $[0,2\pi)$, than $f'$
also has at least $m$ zeroes in this interval (just consider $f$
as a function defined on a circle). An advantage of this statement
is that the differentiation can be iterated as many times as we
please without `losing' zeroes, to obtain

\begin{lemma}
\label{Rolle} Let $f:\Real \rightarrow\Real$ be a $C^{\ell}$
($\ell\geq 1$) $2\pi$-periodic function. If $f$ has at least $m$
zeroes in the interval $[0,2\pi)$, then so do $f^{(j)}$, $1\leq
j\leq \ell$.
\end{lemma}

We now introduce some useful notation. Given a continuous
$2\pi$-periodic function $f:\Real\rightarrow\Real$ with mean-value
$0$:
$$<f>=\frac{1}{2\pi}\int_0^{2\pi}{f(x)dx}=0,$$
we denote by $f^{(-1)}$ the {\it{antiderivative}} of $f$. We fix
the constant of integration by assuming that $<f^{(-1)}>=0$. We
note that the assumption that $<f>=0$ implies that $f^{(-1)}$ is
also $2\pi$-periodic.  For each negative integer $-\ell$ we
define, recursively, $f^{(-\ell)}=(f^{(-\ell+1)})^{(-1)}$, so that
$f^{(-\ell)}$ is the $\ell$-th antiderivative of $f$.

The strategy of our proof of the Sturm-Hurwitz theorem is as
follows: we assume that $n\geq 1$ (otherwise the result is
trivial), which means that $<f>=0$. We will show that for
{\it{sufficiently large}} $\ell$, the $2\pi$-periodic function
$f^{(-\ell)}$ has at least $2n$ zeroes. By lemma \ref{Rolle}, this
implies that $f=(f^{(-\ell)})^{(\ell)}$ has at least $2n$ zeroes,
which is what we want to show.

To prove the above claim we note that when $\ell$ is a multiple of
$4$, the Fourier series of $f^{(-\ell)}$ is
$$f^{(-\ell)}(x)= \sum_{k=n}^{\infty}{\frac{1}{k^\ell}[a_k \cos(kx) +
b_k \sin(kx)]}$$ (when $\ell$ is not a multiple of $4$ we can also
write down the Fourier series of $f^{(-\ell)}$, but for our
purposes we can just assume henceforth that $\ell$ is a multiple
of $4$). We note that the series on the right-hand side indeed
uniformly converges to the function on the left-hand side, since
that function is continuously differentiable. A key observation
for our proof is that when $\ell$ is large, the first term in the
above series becomes `dominant' in such a way as to force the
function $f^{(-\ell)}$ to have $2n$ zeroes. To explain what we
mean by this, let us define
$$g(x)=a_n \cos(nx)+b_n \sin(nx).$$
Since $g^{(-\ell)}(x)=\frac{1}{n^\ell}[a_n \cos(nx)+b_n \sin(nx)]$
is simply a translate of the function
$\frac{1}{n^\ell}\sqrt{a_n^2+b_n^2}\sin(nx)$, the maximum and the
minimum values of $g^{(-\ell)}(x)$ are $\pm
\frac{1}{n^\ell}\sqrt{a_n^2+b_n^2}$, and each is attained $n$
times in the interval $[0,2\pi)$, with maxima and minima
alternating. Thus if we can show that, for sufficiently large
$\ell$, we have
\begin{equation}
\label{ineq}d_\ell\equiv \max_{x\in
\Real}{|f^{(-\ell)}(x)-g^{(-\ell)}(x)|}<\frac{1}{n^\ell}\sqrt{a_n^2+b_n^2},
\end{equation}
then we conclude that $f^{(-\ell)}$ is positive at the maxima of
$g^{(-\ell)}$ and negative at its minima, so that by the
intermediate value theorem $f^{(-\ell)}$ vanishes $2n$ times in
$[0,2\pi)$.

To prove that (\ref{ineq}) holds for sufficiently large $\ell$, we
set $M=\max_{x\in\Real}{|f(x)|}$ and we note that from the formula
for Fourier coefficients we have that $|a_k|,|b_k|<2M$ for all
$k$, so that we can estimate $d_\ell$ from above as follows:
\begin{eqnarray}
\label{ii}d_\ell&=&\max_{x\in
\Real}{\Big{|}\sum_{k=n+1}^{\infty}{\frac{1}{k^\ell}[a_k \cos(kx)
+ b_k \sin(kx)]}\Big{|}}\leq
4M\sum_{k=n+1}^{\infty}{\frac{1}{k^\ell}}\\&\leq&
4M\int_n^{\infty}{\frac{du}{u^{\ell}}}=\frac{4M}{(\ell-1)n^{\ell-1}}.
\nonumber
\end{eqnarray}
It is easy to check that the right-hand side of (\ref{ii}) is
smaller than the right hand side of (\ref{ineq}) for sufficiently
large $\ell$, implying that (\ref{ineq}) holds and completing our
proof.

\vspace{0.5cm} Finally, we remark on some analogues of the
Sturm-Hurwitz theorem, for non-periodic functions. One direction
is to consider trigonometric sums of the form
\begin{equation}
\label{qp} f(x)=\sum_{i=1}^N[a_i \cos(\lambda_i x)+ b_i
\sin(\lambda_i x)],
\end{equation}
where $0<\lambda_1<\lambda_2<...<\lambda_N$ or, more generally,
almost-periodic functions. Analogues of the Sturm-Hurwitz theorem
give a lower bound on the `density' of zeroes in terms of
$\lambda_1$. Such results are closely related to the theory of
`mean motions' \cite{Jessen}. The method used in our proof of the
Sturm-Hurwitz theorem can indeed be adapted to prove such a result
for functions of the form (\ref{qp}), and the interested reader
can do this as an exercise.

Another direction is to assume that the Fourier transform of a
function $f:\Real\rightarrow\Real$ vanishes in an interval
$(-a,a)$ and deduce results about the density of zeroes. This is
investigated in the recent paper \cite{en}. It is interesting to
note that the authors of \cite{en} use some of the ideas used in
previous proofs of the Sturm-Hurwitz theorem, enumerated at the
beginning of this note. It would be interesting to know whether
the method of proof used here can be adapted to yield results in
this direction.


\begin{thebibliography}{9}



\bibitem{Arnold2} V.I. Arnold,
{\it{Topological problems in wave propagation theory and
topological economy principle in algebraic geometry}}, In:
Proceedings in honour to V.I. Arnold for his 60th birthday, AMS
Fields Inst. Com. {\bf{24}}, 39-54 (1999).

\bibitem{Arnold1}  V.I. Arnold,
{\it{Astroidal geometry of hypocycloids and the Hessian topology
of hyperbolic polynomials}}, Russian Math. Surveys {\bf{56}}
(2001), 1019-1083.

\bibitem{en} A. Eremenko \& D. Novikov, Oscillation of Fourier
integrals with a spectral gap, J. Math. Pure Appl., to appear.

\bibitem{Hurwitz} A. Hurwitz, {\it{\"Uber die Fourierschen konstanten
integrierbarer funktionen}}, Math. Ann. {\bf{57}} (1903), 425-446.

\bibitem{Jessen} B. Jessen \& H. Tornehave, {\it{Mean motions and zeroes
of almost periodic functions}}, Acta Math {\bf{77}} (1945),
137-279.

\bibitem{Maure} Y. Martinez-Maure,
{\it{Les multih\'erissons et le th\'eor\`eme de Sturm-Hurwitz}},
Arch. Math. {\bf{80}} (2003), 79-86.

\bibitem{Polya} G. P\'olya, {\it{Qualitatives \"uber W\"armeausgleich}}, Z.
angew. Math. Mech. {\bf{13}} (1933), 125-128.

\bibitem{ps} G. P\'olya \& G. Szeg\"o, `Problems and Theorems in
Analysis', Springer (Berlin), 1972.

\bibitem{Sturm} C. Sturm, {\it{Sur une classe d'equations \`a
diff\'erences partielles}}, J. Math. Pure Appl. {\bf{1}} (1836),
373-444.



\end{thebibliography}
\end{document}